\documentclass[12pt]{article}

\usepackage{amstext    }
\usepackage{amsthm    }
\usepackage{a4}
\usepackage[mathscr]{eucal}
\usepackage{mathrsfs}

\usepackage{amsmath}
\usepackage{amssymb}
\usepackage{amscd}

\newtheorem{theorem}{Theorem}[section]

\newtheorem{proposition}[theorem]{Proposition}
\newtheorem{corollary}[theorem]{Corollary}
\newtheorem{lemma}[theorem]{Lemma}
\newtheorem{remark}[theorem]{Remark}

\newcommand{\cali}[1]{\mathscr{#1}}

\newcommand{\supp}{{\rm supp}}

\newcommand{\loc}{{loc}}
\newcommand{\ddc}{dd^c}
\newcommand{\dc}{d^c}

\newcommand{\DSH}{{\rm DSH}}

\newcommand{\Bc}{\cali{B}}
\newcommand{\Cc}{\cali{C}}

\newcommand{\Fc}{\cali{F}}
\newcommand{\Gc}{\cali{G}}

\newcommand{\Nc}{\cali{N}}

\newcommand{\Et}{{\tt E}}
\newcommand{\St}{{\tt S}}

\newcommand{\FS}{{\rm FS}}

\newcommand{\C}{\mathbb{C}}

\newcommand{\R}{\mathbb{R}}
\renewcommand{\P}{\mathbb{P}}

\title{Exponential estimates for plurisubharmonic functions and 
  stochastic dynamics}

\author{Tien-Cuong Dinh, Vi{\^e}t-Anh Nguy{\^e}n and Nessim Sibony}

\begin{document}

\maketitle

\begin{abstract}
We prove  exponential estimates for plurisubharmonic functions with
respect to Monge-Amp{\`e}re measures with H{\"o}lder continuous
potential. As an application, we
obtain several stochastic properties  for the equilibrium measures associated
to  holomorphic maps on projective spaces. More precisely, we prove 
the exponential decay of correlations, the central limit
theorem for general d.s.h. observables, and the large deviations theorem
for bounded d.s.h. observables and H{\"o}lder continuous observables.
\end{abstract}

\noindent
{\bf Classification AMS 2000:} 32U, 37F.


\section{Introduction}
 
In this paper we prove exponential estimates for plurisubharmonic functions with respect
to a class of probability measures which contains the measures of maximal entropy for many
dynamical systems in several complex variables. This permits to prove
the large deviations
theorem for these dynamical systems and also sharp decay of
correlation estimates. The
results seem to be new even in dimension one.
This type of exponential estimates should play a role in the study of stochastic
properties of dynamical systems in the complex domain.

Let $X$ be a complex manifold of dimension $k$ and $K$ a compact
subset of $X$. Let $\mu$ be a positive measure on $X$. If $\psi$ is a plurisubharmonic function and
if $\mu$ is given by a 
differential form with coefficients in $L^p_\loc$, $p>1$, then  $e^{-\alpha\psi}$ restricted to $K$ is
integrable with respect to $\mu$ for some constant $\alpha>0$. The case where $X$ is an open
set in $\C^k$ and $\mu$ is the Lebesgue measure is a classical
result, see H{\"o}rmander \cite{Hormander} and Skoda \cite{Skoda}. The
general case is a direct consequence. These estimates are
very useful in complex geometry, see e.g. Demailly's book
\cite{Demailly} and the references therein.
They are also very useful in K{\"a}hler-Einstein
geometry and have been developed by Tian-Yau \cite{Tian, TianYau, Yau}.

In this paper, we consider a class of measures satisfying an analogous
property. We first recall some notions, see \cite{DinhSibony1}. The measure $\mu$ is said
to be {\it  locally moderate} if for any open set $U\subset X$,
any compact set $K\subset U$ and any compact familly $\Fc$ of
plurisubharmonic functions (p.s.h. for short) on $U$, there are
constants $\alpha>0$ and $c>0$ such that 
$$\int_Ke^{-\alpha\psi}d\mu\leq c\quad \mbox{for } \psi\in\Fc.$$
This inequality implies that $\Fc$ is bounded in $L_\loc^p(\mu)$ for
$1\leq p<\infty$. In particular, $\mu$ has no mass on pluripolar sets.
The existence of $c$ and $\alpha$ is equivalent to the
existence of $c'>0$ and $\alpha'>0$ satisfying
$$\mu\{z\in K,\ \psi(z)<-M\}\leq c'e^{-\alpha'M}$$
for $M\geq 0$ and $\psi\in\Fc$. Note that the functions on $\Fc$ are
uniformly bounded from above on $K$, see e.g. \cite{Demailly}.
Applying the above estimates to $\log\|z-a\|$, we obtain that the $\mu$-measure
of a ball of center $a\in K$ and of small radius $r$ is bounded by
$r^{\alpha''}$ for some $\alpha''>0$. In the one variable case, this
property is equivalent to the fact that $\mu$ is locally moderate.

Fix a hermitian form $\omega$, i.e. a smooth strictly positive $(1,1)$-form,
on $X$. Let $S$ be a positive closed current of bidegree
$(p,p)$ on $X$. Define the {\it trace measure} of $S$ by $\sigma_S:=S\wedge
\omega^{k-p}$. We say that $S$ is 
{\it locally moderate} if its trace measure 
is locally moderate. So, if $S$ is given by a continuous differential form then
it is locally moderate. Observe that the notion of locally moderate
current does not depend on the choice of $\omega$.

Consider a continuous real-valued function $u$ on the support
$\supp(S)$ of $S$. The multiplication $uS$ defines a current on $X$,
so the current $\ddc(uS)$ is also well-defined. The function $u$ is {\it
  $S$-p.s.h.} if $\ddc(uS)$ is a positive current. 
If $R$ is a positive closed $(1,1)$-current on $X$, we can
locally  write $R=\ddc u$ where $u$ is a p.s.h. function. We call $u$ {\it
 a local potential} of $R$. If $R$ has local continuous potentials then
the wedge-product $R\wedge S$ is well-defined and is locally given by
$R\wedge S:=\ddc (uS)$. Indeed, it is enough to have that $u$ is
locally integrable with respect to the trace measure of $S$.
The wedge-product is a positive closed $(p+1,p+1)$-current
which does not depend on the choice of $u$.
We refer the reader
to \cite{BedfordTaylor, Demailly, FornaessSibony2, BerndtssonSibony, DinhSibony8} for
the intersection theory of currents. Here is one of our main result.

\begin{theorem} \label{th_moderate}
Let $S$ be a locally moderate positive closed
  $(p,p)$-current on a complex manifold $X$. 
If $u$ is a H{\"o}lder continuous $S$-p.s.h. function, then $\ddc(uS)$
is locally moderate. In particular, if $R$ is a positive closed
$(1,1)$-current with H{\"o}lder continuous local potentials, then
$R\wedge S$ is locally moderate.
\end{theorem}

If $u$ is a continuous p.s.h. function on $X$, the Monge-Amp{\`e}re
$(p,p)$-currents, $1\leq p\leq k$, associated to $u$ is defined by induction
$$(\ddc u)^p:=\ddc u\wedge \ldots\wedge \ddc u\quad (p\mbox{
  times}).$$
These currents are very useful in complex analysis and geometry. We have
the following corollary.

\begin{corollary}
Let $u$ be a H{\"o}lder continuous p.s.h. function on $X$. Then the
Monge-Amp{\`e}re currents $(\ddc u)^p$ are locally moderate.
\end{corollary}

We give now
an application to dynamics.
Consider a non-invertible holomorphic endomorphism $f$ of the
projective space $\P^k$. Let $d\geq 2$ denote the algebraic degree of
$f$. That is, $f$ is induced by a homogeneous polynomial endomorphism
of degree $d$ on $\C^{k+1}$. If $V$ is a subvariety of pure
codimension $p$ on $\P^k$, then $f^{-1}(V)$ is a subvariety of pure
codimension $p$ and of degree (counted with multiplicity)
$d^p\deg(V)$. More generally, if $S$ is a positive closed
$(p,p)$-current on $\P^k$ then $f^*(S)$ is a well-defined positive
closed $(p,p)$-current of mass $d^p\|S\|$. Here, we consider the
metric on $\P^k$ induced by the Fubini-Study form $\omega_\FS$ that we
normalize by $\int_{\P^k}\omega_\FS^k=1$. The mass of $S$ is given by
$\|S\|:=\langle S,\omega_\FS^{k-p}\rangle$. We refer the reader to
\cite{Meo, DinhSibony6} for the definition of the pull-back operator
$f^*$ on positive closed currents.

Recall some dynamical properties of $f$ and its iterates
$f^n:=f\circ\cdots \circ f$, $n$ times, see e.g. the survey article
\cite{Sibony}. One can associate to $f$ some canonical invariant
currents. Indeed, $d^{-n} (f^n)^*(\omega_\FS)$ converge to  a
positive closed $(1,1)$-current $T$ of mass 1 on $\P^k$.
The current $T$ has locally H{\"o}lder continuous potentials. So, one can
define $T^p:=T\wedge\ldots\wedge T$, $p$ times. The
currents $T^p$ are {\it the Green currents} associated to $f$ and
$\mu:=T^k$ is the {\it Green measure} of $f$. They are {\it totally
  invariant} by $f$: $d^{-p}f^*(T^p)=d^{-k+p}f_*(T^p)= T^p$. 

\begin{corollary} \label{corollary_dynamics}
Let $f$ be a non-invertible holomorphic endomorphism
  of $\P^k$. Then the Green currents and the Green measure associated
  to $f$ are locally moderate.
\end{corollary} 

This property of the Green measure $\mu$ allows us to prove the
central limit theorem and the large deviations theorem for a large class of observables. Recall that a
{\it quasi-p.s.h. function} on $\P^k$ is locally the difference of a
p.s.h. function and a smooth function. A function is called {\it d.s.h.}
if it is equal outside a pluripolar set to the difference of two
quasi-p.s.h. functions. We identify two d.s.h. functions if
they are equal out of a pluripolar set. Moreover, d.s.h. functions are integrable with respect
to $\mu$, see e.g. \cite{Demailly, DinhSibony3} and Section \ref{section_mixing}. Recall that $\mu$ 
has no mass on pluripolar sets. 

\begin{corollary} \label{corollary_clt}
Let $f$ be a non-invertible holomorphic endomorphism of $\P^k$. Let
$\mu$ denote its Green measure. If a d.s.h. function $\psi$ on
$\P^k$ satisfies $\langle\mu,\psi\rangle=0$ and is not a coboundary,
then it satisfies the central limit theorem with respect to $\mu$.
\end{corollary}

The reader will find more details in Sections \ref{section_mixing} and
\ref{section_clt}. Corollary \ref{corollary_clt}
was known for $\psi$ bounded d.s.h., and for
$\psi$ H{\"o}lder continuous, see \cite{FornaessSibony1, DinhSibony1,
  DinhSibony3, DinhSibony5}. The result was recently
extended by Dupont to unbounded functions $\psi$ with analytic
singularities such that $e^\psi$ is H{\"o}lder continuous \cite{Dupont}.
He uses Ibragimov's approach and gives an application of the
central limit theorem, see also \cite{Dupont2}. Our result relies on the verification of
Gordin's condition in Theorem \ref{th_gordin}. 
Some finer stochastic properties of
$\mu$ (the almost-sure invariance principle, the Donsker and Strassen
principles and the law of the iterated logarithm) can be deduced from
the so-called Philipp-Stout's condition proved by Dupont \cite{Dupont} or
Gordin's condition that we obtain here, see \cite{HallHeyde}. We refer
to \cite{Denker, Haydn, Makarov, DinhNguyenSibony1, PrzytyckiRivera, XiaFu} and the references
therein for some results in the case of dimension 1.
We will prove in Section \ref{section_ldt} that bounded
d.s.h. functions and H{\"o}lder continuous functions 
satisfy the large deviations theorem.

Note that Corollary \ref{corollary_dynamics} 
can be extended to several situations, in particular to H{\'e}non
maps, to regular polynomial automorphisms and also to automorphisms of compact
K{\"a}hler manifolds \cite{Sibony, DinhSibony4}. For a simple proof of the
H{\"o}lder continuity of Green functions, see \cite[Lemma 5.4.2]{DinhSibony8}.

\section{Locally moderate currents}

In this section, we give the proof of Theorem
\ref{th_moderate}. Assume that $S$ is moderate and $u$ is a H{\"o}lder continuous
function on $\supp(S)$ with H{\"o}lder exponent $0<\nu\leq 1$. 

The problem is local. So, we can assume that $U=X$ is the ball $B_2$ of center
0 and of radius 2 in $\C^k$, $K$ is the closed ball $\overline B_{1/2}$ of radius $1/2$
and $\omega$ is the canonical K{\"a}hler
form $\ddc\|z\|^2$. Here, $z=(z_1,\ldots,z_k)$ is a coordinate system
of $\C^k$. We replace $S$ by $S\wedge \omega^{k-p-1}$ in order to
assume that $S$ is of bidegree $(k-1,k-1)$. 
We have the following lemma.

\begin{lemma} \label{lemma_int}
Let $\Gc$ be a compact family of p.s.h. functions on $B_2$. Then $\Gc$
is bounded in $L^1_\loc(\sigma_S)$. 
Moreover, the mass of the measure $\ddc\varphi\wedge S$ is
locally bounded on $B_2$, uniformly on $\varphi\in\Gc$.
\end{lemma}
\proof
Observe that on any compact set $H$ of $B_2$, the functions of $\Gc$ are bounded from above by
the same constant. Subtracting from these functions a constant allows to assume
that they are negative on $H$. We deduce from the 
fact that $S$ is locally moderate that $\int_H \varphi d\sigma_S$ is bounded
uniformly on $\varphi\in\Gc$. 
Indeed, there is a constant $c_H>0$ such that $|\int_H\varphi
d\sigma_S|\leq c_H\|\varphi\|_{L^1(B_2)}$ for negative
p.s.h. functions $\varphi$ on $B_2$. 
This proves the first assertion. 

For the second assertion, consider a compact set $H\subset B_2$. Let
$0\leq\chi\leq 1$ be a cut-off function, smooth, supported on a compact set
$L\subset B_2$ and equal to 1
on $H$. The mass of
$\ddc\varphi\wedge S$ is bounded by the following integral
$$\int\chi \ddc(\varphi S)=\int \ddc\chi \wedge \varphi S\leq
\|\chi\|_{\Cc^2}\int_L|\varphi|d\sigma_S.$$
We have seen that the last term is bounded uniformly on $\varphi$. 
\endproof

We will use the following classical lemma.

\begin{lemma} Let $u$ be a $\nu$-H{\"o}lder continuous function on a
  closed subset $F$ of $B_2$. Then
$u$ can be extended to a $\nu$-H{\"o}lder continuous function on $B_1$.
\end{lemma}
\proof
We can assume that $|u|\leq 1$. Define for $x\in B_1$
$$\widetilde u(x):=\min\{u(y)+A\|x-y\|^\nu,\ y\in F\cap
\overline B_1\}$$
where $A>0$ is  a constant large enough so that 
$|u(a)-u(b)|\leq A\|a-b\|^\nu$ on $F\cap \overline B_1$.
It follows that $\widetilde u(x)=u(x)$ for $x\in F\cap B_1$. We
only have to check that $\widetilde u$ is $\nu$-H{\"o}lder continuous.

Consider two points $x$ and $x'$ in $B_1$ such that $\widetilde u(x)\leq \widetilde
u(x')$. The aim is to bound 
$\widetilde u(x')-\widetilde u(x)$. Let
$y$ be a point in $F\cap \overline B_1$ such that $\widetilde
u(x)=u(y)+A\|x-y\|^\nu$.
By definition of $\widetilde u(x')$, we have
\begin{eqnarray*}
\widetilde u(x')-\widetilde u(x) & \leq &
u(y)+A\|x'-y\|^\nu-\widetilde u(x) \\
& \leq & A\|x'-y\|^\nu -A\|x-y\|^\nu\\
& \leq & A(\|x'-x\|+\|x-y\|)^\nu -A\|x-y\|^\nu\\
& \leq & A\|x'-x\|^\nu
\end{eqnarray*}
since $t\mapsto t^\nu$ is concave increasing on $t\in \R^+$ for $0<\nu\leq
1$. This completes the proof.
\endproof

We continue the proof of Theorem \ref{th_moderate}. 
For simplicity, let $u$ denote the extension of $u$ to $B_1$ as
above for $F=\supp(S)$. Subtracting from $u$ a constant allows to assume that $u\leq
-1$ on $B_1$. 
Define $v(z):=\max(u(z),A\log\|z\|)$ for a constant $A$ large
enough. Observe that since $A$ is large, $v$ is equal to $u$ on
$B_{2/3}$ and to $A\log\|z\|$ near the boundary of $B_1$. 
Moreover, $v$ is $\nu$-H{\"o}lder continuous.
We are interested in an estimate on $\overline B_{1/2}$. So, replacing $u$ by $v$
allows us to assume that $u=A\log\|z\|$ on $B_1\setminus B_{1-4r}$  for some
constant $0<r<1/16$. Fix a smooth function $\chi$
with compact support in $B_{1-r}$, equal to 1 on $B_{1-2r}$ and such that
$0\leq\chi\leq 1$. 

\begin{lemma} \label{lemma_by_parts}
If $\varphi$ is a p.s.h. function on $B_2$ then
\begin{eqnarray*} 
\int_{B_1}\chi\varphi\ddc(u S) & = & -\int_{B_{1-r}\setminus
  B_{1-3r}}\ddc\chi\wedge \varphi uS -
\int_{B_{1-r}\setminus B_{1-3r}} d\chi\wedge \varphi \dc u\wedge S \\
& & +\int_{B_{1-r}\setminus B_{1-3r}}
\dc\chi\wedge \varphi du \wedge S  +\int_{B_{1-r}} \chi u\ddc\varphi \wedge S.
\end{eqnarray*}
\end{lemma}
\proof
Observe that $u$ is smooth on $B_1\setminus B_{1-4r}$. So, all the
previous integrals make sense, see also Lemma \ref{lemma_int}. On the
other hand, one can approximate $\varphi$ by a deacreasing sequence of
smooth p.s.h. functions (one reduce slightly $B_2$ if
necessary). Therefore, it is enough to prove the lemma for $\varphi$
smooth. A direct computation gives
\begin{eqnarray*} 
\int_{B_1}\chi\varphi\ddc(u S) & = & -\int_{B_1}\ddc\chi\wedge \varphi uS-
\int_{B_1} d\chi\wedge \varphi \dc (u S) \\
& & +\int_{B_1}
\dc\chi\wedge \varphi d(u S)  +\int_{B_1} \chi \ddc\varphi \wedge uS.
\end{eqnarray*}
The fact that $\ddc\chi$, $d\chi$, $\dc\chi$ are supported in
$B_{1-r}\setminus B_{1-2r}$ and $\chi$ is supported in $B_{1-r}$ imply the result.
\endproof

\noindent
{\bf End of the proof of Theorem \ref{th_moderate}.} Since $\Fc$ is
compact, it is locally bounded from above.
Subtracting from each function $\varphi\in\Fc$ a constant allows to assume that
$\varphi\leq 0$ on $B_1$.
Define $\varphi_M:=\max(\varphi,-M)$ and
$\psi_M:=\varphi_{M-1}-\varphi_M$ for $\varphi\in \Fc$ and $M\geq
0$. Let $\Gc$ denote the family of all these functions
$\varphi_M$. This family is compact in $L^1_\loc$.
Lemma \ref{lemma_int} implies that in $B_1$ the masses of
$\ddc\varphi_M$ and of $\ddc \varphi_M\wedge S$ are locally bounded
independently of $\varphi\in\Fc$ and of $M$. 
The function $\psi_M$ is positive, bounded by 1, supported in $\{\varphi<-M+1\}$,
and equal to 1 on $\{\varphi<-M\}$. The mass of $\ddc(uS)$ on $\{\varphi<-M\}$ is bounded
by 
$$\int\chi\psi_M \ddc (uS)$$ 
with $\chi$ as in Lemma \ref{lemma_by_parts} above.
We will show that this integral is $\lesssim e^{-\alpha\nu M/3}$ for
some $\alpha>0$. This
implies the result. 

Since $S$ is locally moderate, we have the following estimate for some
$\alpha>0$
\begin{equation*}
\sigma_S\big\{z\in B_{1-r},\ \varphi(z)<-M+1\big\}\lesssim e^{-\alpha M}.
\end{equation*}
Lemma \ref{lemma_by_parts} implies that
\begin{eqnarray*} 
\int_{B_1}\chi\psi_M\ddc(u S) & = & -\int_{B_{1-r}\setminus
  B_{1-3r}}\ddc\chi\wedge \psi_M uS-
\int_{B_{1-r}\setminus B_{1-3r}} d\chi\wedge \psi_M \dc u\wedge S \\
& & +\int_{B_{1-r}\setminus B_{1-3r}}
\dc\chi\wedge \psi_M du \wedge S  +\int_{B_{1-r}} \chi u\ddc\psi_M \wedge S.
\end{eqnarray*}
The first three integrals on
the right hand side are $\lesssim
e^{-\alpha M}$. This is a consequence of the above estimate on $\sigma_S$
and the smoothness of $u$ on $B_1\setminus B_{1-4r}$.
It remains to estimate the last integral.

We use now the $\nu$-H{\"o}lder continuity of $u$.
Define $\epsilon:=e^{-\alpha M/3}$. This is a small constant since we
only have to consider $M$ big.
Write $u=u_\epsilon + (u-u_\epsilon)$ where $u_\epsilon$ is defined on
$B_{1-r}$ and is obtained
from $u$ by convolution with a smooth approximation of identity. The convolution can be chosen so that
$\|u_\epsilon\|_{\Cc^2}\lesssim \epsilon^{-2}$ and
$\|u-u_\epsilon\|_\infty\lesssim \epsilon^\nu=e^{-\alpha\nu M/3}$. Moreover, the
$\Cc^2$-norm of $u_\epsilon$ on $B_{1-r}\setminus B_{1-3r}$ is bounded
independently of $\epsilon$ since $u=\log\|z\|$ on $B_1\setminus
B_{1-4r}$. We have
\begin{eqnarray*}
\int \chi u \ddc \psi_M \wedge S   & = &   \int \chi \ddc \psi_M\wedge  Su_\epsilon +
\int \chi \ddc\psi_M\wedge S (u-u_\epsilon)\\
 & = &   \int \chi \ddc \psi_M\wedge  Su_\epsilon +
\int \chi (\ddc\varphi_{M-1}-\ddc\varphi_M)\wedge S (u-u_\epsilon).
\end{eqnarray*}
By Lemma \ref{lemma_int}, the last integral is $\lesssim
\|u-u_\epsilon\|_\infty\lesssim e^{-\alpha\nu M/3}$. 

Using an expansion as above, we obtain
\begin{eqnarray*}
\int \chi \ddc \psi_M\wedge  Su_\epsilon 
& = & 
\int_{B_{1-r}\setminus B_{1-3r}} \ddc \chi \wedge \psi_M
S u_\epsilon + \int_{B_{1-r}\setminus B_{1-3r}} d \chi \wedge \psi_M
S\wedge \dc u_\epsilon \\
&& -\int_{B_{1-r}\setminus B_{1-3r}} \dc \chi \wedge \psi_M
S \wedge du_\epsilon + \int \chi \psi_M  S \wedge \ddc u_\epsilon.
\end{eqnarray*}
As above, the first three integrals on the right hand side are $\lesssim
e^{-\alpha M}$ because $u_\epsilon$ has bounded $\Cc^2$-norm on $B_{1-r}\setminus
B_{1-3r}$. Consider the last integral. Since $\psi_M$ is supported in
$\{\varphi\leq -M+1\}$, the estimate on $\sigma_S$ implies that the
considered integral is
$\lesssim e^{-\alpha M}\|u_\epsilon\|_{\Cc^2}\lesssim e^{- \alpha
  M}\epsilon^{-2}= e^{-\alpha\nu M/3}$. We deduce from all
the previous estimates that
$$\int \chi u\ddc \psi_M \wedge S \lesssim e^{-\alpha\nu M/3}.$$
This completes the proof.
\hfill $\square$

\begin{remark} \rm
On a compact K{\"a}hler manifold $X$, one can introduce the notion of
{\it (globally) moderate current}. For this purpose, in the
definition, one replaces 
local p.s.h. functions by (global) quasi-p.s.h. functions. In the case
where $X=\P^k$, the first and third authors introduced in \cite{DinhSibony8} a notion of
super-potential for positive closed $(p,p)$-currents. One can prove
that currents with H{\"o}lder continuous super-potentials are moderate
and the intersection of currents with H{\"o}lder continuous
super-potentials admits H{\"o}lder continuous super-potentials. If a
$(p,p)$-current admits a H{\"o}lder continuous potential, it has
a H{\"o}lder continuous super-potential and then is moderate. 
\end{remark}

\section{Decay of correlations} \label{section_mixing}

Let $\mu$ be the
Green measure of an endomorphism $f$ of algebraic degree $d\geq 2$ of $\P^k$.
In this section, we will prove that $\mu$ is {\it mixing} and {\it
  exponentially mixing} in different senses.
If $\phi$ is a
d.s.h. function on $\P^k$ we can write $\ddc\phi=R^+-R^-$ where
$R^\pm$ are positive closed $(1,1)$-currents. The {\it d.s.h. norm} of
$\phi$ is defined by
$$\|\phi\|_\DSH:=\|\phi\|_{L^1(\P^k)}+\inf\|R^\pm\|$$
with $R^\pm$ as above. Note that $R^+$ and $R^-$ have the same mass
since they are cohomologous, and that $\|\cdot\|_\DSH\lesssim \|\cdot\|_{\Cc^2}$. 
The following result 
was proved in \cite{FornaessSibony1, DinhSibony1,  DinhSibony3} for $p=+\infty$.

\begin{theorem} \label{th_mixing}
Let $f$ be a holomorphic endomorphism of algebraic degree $d\geq 2$
and $\mu$ its Green measure. Then for every $1< p\leq +\infty$
there is a constant $c>0$ such that
$$|\langle \mu, (\varphi\circ f^n)\psi\rangle
-\langle\mu,\varphi\rangle \langle\mu,\psi\rangle|\leq c
d^{-n}\|\varphi\|_{L^p(\mu)}\|\psi\|_\DSH$$
for $n\geq 0$, $\varphi$ in $L^p(\mu)$ and $\psi$
d.s.h. Moreover, for $0\leq \nu\leq 2$ there is a constant $c>0$
such that 
$$|\langle \mu, (\varphi\circ f^n)\psi\rangle
-\langle\mu,\varphi\rangle \langle\mu,\psi\rangle|\leq c
d^{-n\nu/2}\|\varphi\|_{L^p(\mu)}\|\psi\|_{\Cc^\nu}$$
for $n\geq 0$, $\varphi$ in $L^p(\mu)$ and $\psi$
of class $\Cc^\nu$.
\end{theorem}

The expression $\langle \mu, (\varphi\circ f^n)\psi\rangle
-\langle\mu,\varphi\rangle \langle\mu,\psi\rangle$ is called {\it the
  correlation of order} $n$ between the observables $\varphi$ and
$\psi$. The measure $\mu$ is said to be {\it mixing} if this correlation
converges to 0 as $n$ tends to infinity, for smooth
observables (or equivalently, for continuous, bounded or $L^2(\mu)$
observables).

Observe that the second assertion in Theorem \ref{th_mixing} is
a consequence of the first one. Indeed, on one hand, since 
$\|\psi\|_\DSH\lesssim\|\psi\|_{\Cc^2}$, we obtain the second
assertion for $\nu=2$. 
On the other hand, we have since $\mu$ is invariant
$$|\langle \mu, (\varphi\circ f^n)\psi\rangle
-\langle\mu,\varphi\rangle \langle\mu,\psi\rangle|\leq 2
\|\varphi\|_{L^1(\mu)}\|\psi\|_{\Cc^0}\lesssim \|\varphi\|_{L^p(\mu)}\|\psi\|_{\Cc^0}. $$
So, the second assertion holds for $\nu=0$.
The theory of interpolation between the Banach spaces $\Cc^0$ and
$\Cc^2$ \cite{Triebel} implies that
$$|\langle \mu, (\varphi\circ f^n)\psi\rangle -\langle\mu,\varphi\rangle \langle\mu,\psi\rangle|\lesssim 
d^{-n\nu/2}\|\varphi\|_{L^p(\mu)}\|\psi\|_{\Cc^\nu}.$$

We prove now the first assertion  in Theorem \ref{th_mixing}. Since $\mu$ is invariant, the
assertion is clear when $\psi$ is constant. 
Indeed, in this case, the correlations vanish.
So, subtracting from
$\psi$ a constant allows to assume that $\langle \mu,\psi\rangle =0$.
We have to bound
$ |\langle \mu, (\varphi\circ f^n)\psi\rangle|$. 
Consider the following {\it weak topology} on the space $\DSH(\P^k)$
of d.s.h. functions. We say that the sequence $(\phi_n)$ converges to
$\phi$ in $\DSH(\P^k)$ if $\phi_n\rightarrow\phi$ in the sense of
currents and if $\|\phi_n\|_\DSH$ is bounded uniformly on $n$. 
Recall here some basic properties of d.s.h. functions, see \cite{DinhSibony3}. 

\begin{proposition} \label{prop_dsh}
Let $\phi$ be a
d.s.h. function. There is a constant $c>0$ independent of $\phi$ and two
negative quasi-p.s.h. functions $\phi^\pm$ such that $\phi=\phi^+-\phi^-$,
$\|\phi^\pm\|_\DSH\leq c\|\phi\|_\DSH$ and $\ddc\phi^\pm\geq
-c\|\phi\|_\DSH\omega_\FS$.  
Moreover, $|\phi|$ is d.s.h. and $\||\phi|\|_\DSH\leq c\|\phi\|_\DSH$.
If $\phi_n\rightarrow\phi$ in $\DSH(\P^k)$ then 
$\phi_n\rightarrow\phi$ in $L^p$ for $1\leq p<+\infty$. 
\end{proposition}

Since $\mu$ is
locally a Monge-Amp{\`e}re measure with continuous potential, $\psi\mapsto
\langle\mu,\psi\rangle$ is continuous with respect to the considered
topology on $\DSH(\P^k)$. We say that $\mu$ is PC. This allows to prove that the
DSH-norm of $\phi$ is equivalent to the following norm
$$\|\phi\|_\DSH':=|\langle\mu,\phi\rangle| + \inf\|R^\pm\|$$
where we write as above $\ddc \phi=R^+-R^-$, see \cite{DinhSibony3}. 
In particular, $\log |h|$ is d.s.h. for any rational function $h$ on
$\P^k$, and similarly 
for the potential $u$ of any positive closed $(1,1)$-current $R$, i.e
a quasi-p.s.h. function $u$ such that $\ddc u=R-c\omega$ for some
constant $c$.

Consider the codimension
1 subspace $\DSH_0(\P^k)$ of $\DSH(\P^k)$ defined by $\langle
\mu,\phi\rangle=0$. On this subspace, one has
$\|\phi\|_{\DSH}'=\inf\|R^\pm\|$. Recall that $\mu$ is totally
  invariant : 
$f^*\mu=d^k\mu$. Then, the space $\DSH_0(\P^k)$ is invariant
under $f_*$. Recall that $f_*\phi$ is defined by
$$f_*\phi(x):=\sum_{y\in f^{-1}(x)} \phi(y)$$
where the points in $f^{-1}(x)$ are counted with multiplicities (there
are exactly $d^k$ points). The mass of a positive closed current on
$\P^k$ can be computed cohomologically. We have
$\|f_*R^\pm\|=d^{k-1}\|R^\pm\|$ and hence $\|f_*\phi\|_\DSH'\leq
d^{k-1}\|\phi\|_\DSH'$ on $\DSH_0(\P^k)$. Define also the {\it Perron-Frobenius operator} by
$$\Lambda\phi:=d^{-k} f_*\phi.$$
Since $\mu$ is totally invariant, this is the adjoint operator of
$f^*$ on $L^2(\mu)$. 
Observe that $\|\Lambda\phi\|'_\DSH\leq d^{-1} \|\phi\|'_\DSH$ on $\DSH_0(\P^k)$. So,
$\Lambda$ has a spectral gap on $\DSH(\P^k)$: the constant
functions correspond to the eigenvalue 1 and the spectral radius on
$\DSH_0(\P^k)$ is bounded by $d^{-1}<1$. 
 
\begin{proposition} \label{prop_iterate_dsh}
There are constants $c>0$ and $\alpha>0$ such that for  $\psi\in\DSH_0(\P^k)$ with
$\|\psi\|_\DSH\leq 1$ and for every $n\geq 0$ we have
$$\langle \mu, e^{\alpha d^n |\Lambda^n\psi|}\rangle \leq c.$$ 
In particular, there is a constant $c>0$
independent of $\psi\in\DSH_0(\P^k)$ such that
$$\|\Lambda^n\psi\|_{L^q(\mu)}\leq cq d^{-n}\|\psi\|_\DSH$$
for every $n\geq 0$ and every  $1\leq q <+\infty$.
\end{proposition}
\proof
Since $\|\cdot\|_\DSH$ and $\|\cdot\|_\DSH'$ are equivalent, we
assume for simplicity that $\|\psi\|'_\DSH=1$. 
Observe that $d^n \Lambda^n\psi$ belongs to the family of
functions in $\DSH_0(\P^k)$ with $\|\cdot\|_\DSH'$ norm bounded by
1. It follows from Proposition \ref{prop_dsh} that  $d^n |\Lambda^n\psi|$ belongs
to a compact  family of d.s.h. functions. By Theorem \ref{th_moderate} and
Proposition \ref{prop_dsh}, there are positive constants $\alpha$ and $c$ such
that 
$$\langle \mu, e^{\alpha d^n |\Lambda^n\psi|}\rangle \leq c.$$ 
Since $e^x\geq x^q/q!$ for $x\geq 0$ and $1\leq q <+\infty$, we
deduce, using the inequality $q!\leq q^q$, 
that $\|d^n\Lambda^n\psi\|_{L^q(\mu)}\leq cq$ for some constant
$c>0$ independent of $\psi$, $n$ and $q$. 
\endproof

\noindent
{\bf End of the proof of Theorem \ref{th_mixing}.}
Let $1< q <+\infty$ such that $p^{-1}+q^{-1}=1$. 
Using a simple coordinate change, Proposition \ref{prop_iterate_dsh} and
the H{\"o}lder inequality, we obtain  that
\begin{eqnarray*}
|\langle \mu,(\varphi\circ f^n)\psi\rangle| & = &  d^{-kn} |\langle (f^n)^*\mu,(\varphi\circ f^n)\psi\rangle|
= |\langle \mu, \varphi\ \Lambda^n\psi\rangle| \\
& \lesssim & \|\varphi\|_{L^p(\mu)} \|\Lambda^n\psi\|_{L^q(\mu)} 
\lesssim d^{-n} \|\varphi\|_{L^p(\mu)}
\|\psi\|_\DSH.
\end{eqnarray*}
This completes the proof.
\hfill $\square$

\

It is shown in \cite{DinhSibony1} that $\mu$ is mixing of any order and is
K-mixing.  More precisely, we have for
every $\psi$ in $L^2(\mu)$
$$\lim_{n\rightarrow\infty} \sup_{\|\varphi\|_{L^2(\mu)}=1}|\langle
\mu, (\varphi\circ f^n)\psi\rangle -\langle
\mu,\varphi\rangle\langle\mu,\psi\rangle|=0.$$
The reader can deduce the K-mixing from Theorem
\ref{th_mixing} and the fact that $\Lambda$ has
norm 1 when it acts on $L^2(\mu)$.

The following result gives the
exponential mixing of any order. It can be extended to 
H{\"o}lder continuous observables using the interpolation theory.

\begin{theorem}
Let $f$, $d$ and $\mu$ be as in Theorem \ref{th_mixing} and $r\geq 1$
an integer. Then there is a constant $c>0$ such that
$$\Big|\langle \mu, \psi_0(\psi_1\circ f^{n_1})\ldots (\psi_r\circ f^{n_r})\rangle
-\prod_{i=0}^r\langle\mu,\psi_i\rangle \Big|\leq c
d^{-n}\prod_{i=0}^r\|\psi_i\|_\DSH$$
for $0=n_0\leq n_1\leq \cdots\leq n_r$, $n:=\min_{0\leq i< r} (n_{i+1}-n_i)$ and $\psi_i$ 
d.s.h. 
\end{theorem}
\proof
The proof is by induction on $r$. The case $r=1$ is a consequence of Theorem
\ref{th_mixing}. Suppose the result is true for $r-1$. We have to
check it for $r$. Without loss of generality, assume that
$\|\psi_i\|_\DSH\leq 1$. This implies that $m:=\langle \mu,\psi_0\rangle$ is
bounded. The invariance of $\mu$ and the hypothesis of induction imply that
\begin{eqnarray*}
\lefteqn{\Big|\langle \mu, m(\psi_1\circ f^{n_1})\ldots (\psi_r\circ f^{n_r})\rangle
-\prod_{i=0}^r\langle\mu,\psi_i\rangle \Big|}\\
&=&  \Big|\langle \mu, m\psi_1(\psi_2\circ f^{n_2-n_1})\ldots (\psi_r\circ f^{n_r-n_1})\rangle
-m\prod_{i=1}^r\langle\mu,\psi_i\rangle \Big|\leq  c d^{-n}
\end{eqnarray*}
for some constant $c>0$. In order to get the desired estimate, it is
enough to show that
$$\Big|\langle \mu, (\psi_0-m)(\psi_1\circ f^{n_1})\ldots 
(\psi_r\circ f^{n_r})\rangle\Big|\leq c d^{-n}.$$
Observe that the operator $(f^n)^*$ acts on $L^p(\mu)$ for $p\geq 1$
and its norm is bounded by 1. Using the invariance of $\mu$ and the
H{\"o}lder inequality, we get for $p:=r+1$
\begin{eqnarray*}
\lefteqn{\Big|\langle \mu, (\psi_0-m)(\psi_1\circ f^{n_1})\ldots 
(\psi_r\circ f^{n_r})\rangle\Big| }\\
&\leq& \Big|\langle \mu, \Lambda^{n_1}(\psi_0-m)\psi_1\ldots 
(\psi_r\circ f^{n_r-n_1})\rangle\Big| \\
&\leq& \|\Lambda^{n_1}(\psi_0-m)\|_{L^p(\mu)}\|\psi_1\|_{L^p(\mu)}\ldots
\|\psi_r\circ f^{n_r-n_1}\|_{L^p(\mu)}\\
& \leq & cd^{-n_1}\|\psi_1\|_{L^p(\mu)}\ldots\|\psi_r\|_{L^p(\mu)},
\end{eqnarray*}
for some constant $c>0$.
Since $\|\psi_i\|_{L^p(\mu)}\lesssim \|\psi_i\|_\DSH$, the previous
estimates imply the result. Note that as in Theorem \ref{th_mixing},
it is enough to assume that $\psi_i$ is d.s.h. for $i\leq r-1$ and
$\psi_r$ is in $L^p(\mu)$ for some $p>1$. 
\endproof

We obtain from Proposition \ref{prop_iterate_dsh} the following result.

\begin{proposition} \label{prop_iterate_holder}
Let $0<\nu\leq 2$ be a constant. There are constants $c>0$ and
$\alpha>0$ such that if $\psi$ is a $\nu$-H{\"o}lder continuous function
with $\|\psi\|_{\Cc^\nu}\leq 1$ and $\langle\mu,\psi\rangle=0$, then 
$$\langle \mu, e^{\alpha d^{n\nu/2} |\Lambda^n\psi|}\rangle \leq c
\qquad \mbox{for every}\qquad n\geq 0.$$ 
Moreover, there is a constant $c>0$
independent of $\psi$  such that
$$\|\Lambda^n\psi\|_{L^q(\mu)}\leq cq^{\nu/2} d^{-n\nu/2}$$
for every $n\geq 0$ and every  $1\leq q <+\infty$.
\end{proposition}
\proof
We only consider the spaces of functions
$\psi$ such that $\langle\mu,\psi\rangle=0$. 
By Proposition \ref{prop_iterate_dsh}, since $\|\cdot\|_\DSH\lesssim
\|\cdot\|_{\Cc^2}$, we have 
$$\|\Lambda^n\psi\|_{L^q(\mu)}\leq cq d^{-n}\|\psi\|_{\Cc^2},$$
with $c>0$ independent of $q$ and of $\psi$. On the other hand, by
definition of $\Lambda$, we have 
$$\|\Lambda^n\psi\|_{L^q(\mu)}\leq
\|\Lambda^n\psi\|_{L^\infty(\mu)}\leq \|\psi\|_{\Cc^0}.$$
The theory of interpolation between the Banach spaces $\Cc^0$ and
$\Cc^2$ \cite{Triebel} (applied to the linear operator
$\psi\mapsto\Lambda^n\psi-\langle\mu,\psi\rangle$) implies that
$$\|\Lambda^n\psi\|_{L^q(\mu)}\leq A_\nu [cq
d^{-n}]^{\nu/2}\|\psi\|_{\Cc^\nu},$$
for some constant $A_\nu>0$ depending only on $\nu$ and on $\P^k$. This
gives the second assertion in the proposition.

For the first assertion, assume that $\|\psi\|_{\Cc^\nu}\leq 1$. Fix a
constant $\alpha>0$ small enough. We have
$$\langle\mu, e^{\alpha d^{n\nu/2} |\Lambda^n\psi|}\rangle  = 
\sum_{q\geq 0} {1\over q!} \langle \mu,
|\alpha d^{n\nu/2}\Lambda^n\psi|^q\rangle 
\leq \sum_{q\geq 0} {1\over q!} \alpha^q c^qq^q.$$
By Stirling's formula, the last sum converges. This implies the result.
\endproof

\section{Central limit theorem} \label{section_clt}

In this section, we give the proof of Corollary
\ref{corollary_clt}. We first recall some facts
\cite{KatokHasselblatt, Walters}. 
Let $(M,\Fc,m)$ be a probability space  and $g:M\rightarrow M$
a measurable map which preserves $m$, i.e. $m$ is $g_*$-invariant :
$g_*m=m$. 
The measure $m$ is {\it ergodic} if for any measurable set $A$ such that
$g^{-1}(A)=A$ we have $m(A)=0$ or $m(A)=1$. This is
equivalent to the property that
$m$ is extremal in the convex set of 
invariant probability measures (if $m$ is mixing then it is
ergodic). When $m$ is ergodic, Birkhoff's theorem  implies that if
$\psi$ is an observable in $L^1(m)$ then 
$$\lim_{n\rightarrow\infty} {1\over n}
\Big[\psi(x)+\psi(g(x))+\cdots+\psi(g^{n-1}(x))\Big]=\langle m,\psi\rangle$$
for $m$-almost every $x$. 

Assume now that
$\langle m,\psi\rangle=0$. Then, the previous limit is equal to 0.
The central limit theorem (CLT for short), when it holds, gives the speed of
this convergence. We say that {\it $\psi$ satisfies the CLT} if
there is a constant $\sigma>0$ such that 
$${1\over \sqrt{n}}
\Big[\psi(x)+\psi(g(x))+\cdots+\psi(g^{n-1}(x))\Big]$$
converges in distribution to the Gaussian random variable $\Nc(0,\sigma)$ of mean 0 and of
variance $\sigma$. Recall that $\psi$ is a {\it coboundary} if there
is a function $\psi'$ in $L^2(\mu)$ such that $\psi=\psi'-\psi'\circ
g$. In this case, one easily checks  that
$$\lim_{n\rightarrow\infty} {1\over \sqrt{n}}
\Big[\psi(x)+\psi(g(x))+\cdots+\psi(g^{n-1}(x))\Big]=\lim_{n\rightarrow\infty} {1\over \sqrt{n}}
\Big[\psi'(x)-\psi'(g^n(x))\Big]= 0$$
in distribution. So, $\psi$ does not satisfies the CLT
(sometimes, one says that it satisfies the CLT for $\sigma=0$). 

The CLT can be deduced from some 
strong mixing, see \cite{BonattiDiazViana, Gordin, Liverani, Viana}.
In the following result, $\Et(\psi|\Fc_n)$ denotes {\it the expectation}
of $\psi$ with respect to $\Fc_n$, that is, $\psi\mapsto \Et(\psi|\Fc_n)$
is the orthogonal projection from $L^2(m)$ onto the subspace generated
by $\Fc_n$-measurable functions.

\begin{theorem}[Gordin] \label{th_gordin}
Consider the decreasing sequence $\Fc_n:=g^{-n}(\Fc)$, $n\geq
0$, of algebras. Let $\psi$ be
a real-valued function in $L^2(m)$ such that
$\langle m,\psi\rangle=0$.  Assume that 
$$\sum_{n\geq 0} \|\Et(\psi |\Fc_n)\|_{L^2(m)}<\infty.$$
Then the  positive number $\sigma$ defined by
$$\sigma^2:=\langle m,\psi^2\rangle + 2\sum_{n\geq 1} \langle
m, \psi(\psi\circ g^n)\rangle$$
is finite. It vanishes if and only if $\psi$ is a coboundary. Moreover, when $\sigma\not=0$,
then $\psi$ satisfies the CLT with variance $\sigma$.  
\end{theorem}

Note that $\sigma$ is equal to the limit of
$n^{-1/2}\|\psi+\cdots+\psi\circ g^{n-1}\|_{L^2(m)}$. The last
expression is equal to $\|\psi\|_{L^2(m)}$ if the family
$(\psi\circ g^n)_{n\geq 0}$ is orthogonal in $L^2(m)$.

We now prove Corollary \ref{corollary_clt}. 
Since $\mu$ is mixing, it is ergodic. So, we
can apply Gordin's theorem to the map $f$ on $(\P^k,\Bc,\mu)$ where $\Bc$
is the canonical Borel algebra.

\begin{lemma} \label{lemma_Borel}
Let $\Bc_n:=f^{-n} (\Bc)$ for $n\geq 0$. Then for $\phi\in L^2(\mu)$
we have
$$\Et(\phi|\Bc_n)=(\Lambda^n\phi)\circ f^n\qquad \mbox{and}\qquad
\|\Et(\phi|\Bc_n)\|_{L^2(\mu)}=\|\Lambda^n\phi\|_{L^2(\mu)}.$$
\end{lemma}
\proof
Consider a function in $L^2(\mu)$ which is measurable with respect to
$\Bc_n$. It has the form $\xi\circ f^n$. We have $\|\xi\circ
f^n\|_{L^2(\mu)}=\|\xi\|_{L^2(\mu)}$ since $\mu$ is invariant. Hence,
$\xi\in L^2(\mu)$.  We deduce from the identity $d^{-kn}
(f^n)^*\mu=\mu$ that
\begin{eqnarray*}
\|\Et(\phi|\Bc_n)\|_{L^2(\mu)} & = & \sup_{\|\xi\|_{L^2(\mu)}=1}
  |\langle\mu,(\xi\circ f^n)\phi\rangle|\\
& = & 
\sup_{\|\xi\|_{L^2(\mu)}=1}
d^{-kn}  |\langle(f^n)^*\mu,(\xi\circ f^n)\phi\rangle|\\
& = & \sup_{\|\xi\|_{L^2(\mu)}=1}
  |\langle \mu,\xi \Lambda^n\phi \rangle| \\ 
& = & \|\Lambda^n\phi\|_{L^2(\mu)}.
\end{eqnarray*}
The computation also shows that the previous supremum 
is reached when $\xi$ is proportional to $\Lambda^n\phi$. It follows that
$\Et(\phi|\Bc_n)=(\Lambda^n\phi)\circ f^n$. 
\endproof

\noindent
{\bf End of the proof of Corollary \ref{corollary_clt}.} By Proposition
\ref{prop_iterate_dsh} and Lemma \ref{lemma_Borel}, since $\psi$ is
d.s.h., we have 
 $\|\Et(\psi|\Bc_n)\|_{L^2(\mu)}\lesssim d^{-n}$. Hence, 
$\sum_{n\geq 0 }\|\Et(\psi|\Bc_n)\|_{L^2(\mu)}$ converges.
It is enough to apply Theorem \ref{th_gordin} in order to get the result. \hfill $\square$

\begin{remark} \rm
If $\psi$ is an observable in $L^\infty(\mu)$, then
$\|\Lambda^n\psi\|_{L^\infty(\mu)} \leq
\|\psi\|_{L^\infty(\mu)}$. Hence, by Lemma \ref{lemma_Borel}, the Gordin's condition in Theorem
\ref{th_gordin} is a consequence of the condition $\sum_{n\geq 1}
\|\Lambda^n\psi\|_{L^1(\mu)}^{1/2}<+\infty$. In particular, H{\"o}lder
continuous observables satisfy the CLT, see Proposition \ref{prop_iterate_holder} and
\cite{DinhSibony3, DinhSibony5} for meromorphic maps. 
\end{remark}

The following proposition gives us the next term in the expansion of
the $L^2$-norm of Birkhoff's sums.

\begin{proposition}
Let $\psi$ be a d.s.h. or an $\nu$-H{\"o}lder continuous
function, with $0<\nu\leq 2$, such that $\langle\mu,\psi\rangle=0$.
Let $\sigma\geq 0$ and $\gamma$ be the constants defined by
$$\sigma^2:=\langle \mu,\psi^2\rangle + 2\sum_{n\geq 1} \langle
\mu, \psi(\psi\circ f^n)\rangle \quad \mbox{and}\quad 
\gamma:=2\sum_{n\geq 1} n\langle\mu, \psi(\psi\circ f^n)\rangle.$$
Then 
$$\|\psi+\cdots+\psi\circ f^{n-1}\|_{L^2(\mu)}^2-n\sigma^2+\gamma$$
is of order $O(d^{-n})$ if $\psi$ is d.s.h. and $O(d^{-n\nu/2})$ if
$\psi$ is $\nu$-H{\"o}lder continuous.
\end{proposition}
\proof
Since $\mu$ is invariant, we have $\langle
\mu, (\psi\circ f^l)(\psi\circ f^m)\rangle =\langle\mu,\psi(\psi\circ
f^{m-l})\rangle$ for $m\geq l$. It follows that
\begin{eqnarray*}
\|\psi+\cdots+\psi\circ f^{n-1}\|_{L^2(\mu)}^2 & = & \sum_{0\leq
  l,m\leq n-1} \langle\mu,(\psi\circ f^l)(\psi\circ f^m)\rangle \\
& = & n\langle \mu,\psi^2\rangle + \sum_{1\leq j\leq n-1} 2(n-j)
\langle \mu,\psi(\psi\circ f^j)\rangle \\
& = & n\sigma^2-\gamma +\sum_{j\geq n}2(j-n) \langle\mu,\psi(\psi\circ f^j)\rangle. 
\end{eqnarray*} 
Theorem \ref{th_mixing} implies the result. Note that this theorem
also implies
that the series in the definition for $\sigma^2$ and for
$\gamma$ are convergent. Moreover, the previous computation gives that
$\sigma^2$ is the limit of 
$n^{-1}\|\psi+\cdots+\psi\circ f^{n-1}\|_{L^2(\mu)}^2$ which is a
positive number. 
\endproof

\section{Large deviations theorem} \label{section_ldt}

In this section, we prove the large deviations theorem (LDT for short)
for the equilibrium measure of holomorphic endomorphisms of $\P^k$. We have
the following result which holds in particular for $\Cc^2$ observables.

\begin{theorem} \label{th_deviation}
Let $f$ be a holomorphic endomorphism of $\P^k$ of algebraic degree
$d\geq 2$. Then the equilibrium measure $\mu$ of $f$ satisfies the
large deviations theorem (LDT) for bounded d.s.h. observables. More
precisely, if $\psi$ is a bounded d.s.h. function then for every
$\epsilon>0$ there is a constant $h_\epsilon>0$ such that 
\begin{equation*}
\mu\Big\{ z\in \P^k:\   \big\vert\frac{1}{n}\sum_{j=0}^{n-1}\psi\circ
f^j(z)- \langle \mu,\psi\rangle \big\vert >\epsilon \Big\}\leq
e^{-n (\log n)^{-2}h_\epsilon}
\end{equation*}
for all $n$ large enough\footnote{in the classical large deviations
  theorem for independent random variables, there is no factor $(\log
  n)^{-2}$ in the previous estimate}.   
\end{theorem}

We start  with the  following Bennett's type inequality, see \cite[Lemma 2.4.1]{DemboZeitouni}.

\begin{lemma}\label{lemma_Bennett}
Let $(M,\Fc,m)$ be  a  probability space and  $\Gc$ a  $\sigma$-subalgebra  of $\Fc.$
Assume that there is a constant $0<\nu<1$ and an element $A\in\Fc$
such that $m(A\cap B)=\nu m(B)$ for every $B\in\Gc$. Define
$s^-:=\max\big\{1, \nu^{-1}(1-\nu)\big\}$ and $s^+:=\max\big\{1,\nu(1-\nu)^{-1}\big\}$. 
Let $\psi$ be a real-valued function on $M$ such that
$\|\psi\|_{L^\infty(m)}\leq b$ and $\Et(\psi|\Gc)=0$.  
Then  
$$\Et(e^{\lambda\psi}|\Gc) \leq \nu e^{-s^-\lambda b}+ (1-\nu)e^{s^+\lambda b}$$ 
for every $\lambda\geq 0$. 
\end{lemma}
\proof  Fix a strictly positive constant $\lambda$. 
Let $\psi_0$ be the function which is equal to
$t^-:=-s^-\lambda b$ on $A$ and to 
$t^+:=s^+\lambda b$ on $M\setminus A$. We have
$\psi_0^2\geq (\lambda b)^2\geq (\lambda \psi)^2$. We deduce from the
hypothesis on $A$ and the relation $-\nu s^-+(1-\nu)s^+=0$ that $\Et(\psi_0|\Gc)=0$. 
Let  $g(t)=a_0t^2+a_1t+a_2$,  be the unique quadratic function such that
$h(t):=g(t)-e^t$ satisfies $h(t^+)=0$ and $h(t^-)=h'(t^-)=0$.
We have $g(\psi_0)=e^{\psi_0}$.

Since $h''(t)=2a_0-e^t$ admits at most one zero, $h'$ admits at most two
zeros. The fact that $h(t^-)=h(t^+)=0$ implies that $h'$ vanishes in
$]t^-,t^+[$. Hence $h'$ admits exactly one zero at $t^-$ and another in
$]t^-,t^+[$. We deduce that $h''$ admits a zero. This implies that
$a_0>0$. Moreover, $h$ vanishes only at $t^-$, $t^+$ and
$h'(t^+)\not=0$. It follows that $h(t)\geq 0$ on $[t^-,t^+]$ because
$h$ is negative near $+\infty$.
Thus, $e^t\leq g(t)$ on $[t^-,t^+]$ and then $\Et(e^{\lambda\psi}|\Gc)\leq  \Et(g(\lambda \psi)|\Gc).$

Since  $a_0>0$, if an observable $\phi$ satisfies $\Et(\phi|\Gc)=0$, 
then $\Et(g(\phi)|\Gc)$ is an increasing function of $\Et(\phi^2|\Gc)$. 
Now, using the properties of $\psi$ and $\psi_0$, 
we obtain 
\begin{eqnarray*}
\Et(e^{\lambda\psi}|\Gc) & \leq &   \Et(g(\lambda \psi)|\Gc)\leq  
\Et(g(\psi_0)|\Gc) \\
 & = & \Et(e^{\psi_0}|\Gc)
 \leq  \nu e^{-s^-\lambda b}+(1-\nu)e^{s^+\lambda b},
\end{eqnarray*} 
which completes the proof.
\endproof

We continue the proof of Theorem \ref{th_deviation}.  
Without loss of generality we may assume  that  $\langle
\mu,\psi\rangle=0$, $|\psi|\leq 1$ and $\|\psi\|_\DSH\leq 1$.
The  general idea is to  write $\psi=\psi'+(\psi''-\psi''\circ f)$
for  functions $\psi'$ and $\psi''$ in  $\DSH_0(\P^k)$  such that 
\begin{equation*}
\Et(\psi'\circ f^n|\Bc_{n+1})=0,\qquad  n\geq 0
\end{equation*}
where  $\Bc_n:=f^{-n}(\Bc)$ as above with $\Bc$ the canonical Borel
algebra of $\P^k$.
In the  language of probability theory, these identities
mean that  $(\psi'\circ f^n)_{n\geq 0}$ is  a {\it reversed martingale
  difference} as in Gordin's approach, see also \cite{Viana}.  
The strategy  is  to  prove the LDT  for  $\psi'$ and
for the coboundary $\psi''-\psi''\circ f$. Theorem \ref{th_deviation} is
in fact a consequence of Lemmas \ref{lemma_ldt} and \ref{lemma_ldt_bis} below.  

Define
\begin{equation*}
\psi'':=-\sum_{n=1}^{\infty} \Lambda^n\psi,\qquad  \psi':=\psi-(\psi''-\psi''\circ f).
\end{equation*}  
Using the estimate  $ \|\Lambda\phi\|'_{\DSH}\leq d^{-1}\|
\phi\|'_{\DSH}$ on $\DSH_0(\P^k)$,  
we  see that  $\psi'$ and $\psi''$ are in $\DSH_0(\P^k)$ with norms
bounded by some constant. In
particular, they belong to $L^2(\mu)$. However, we loose the boundedness:
these functions are not necessarily in $L^\infty(\mu)$.

\begin{lemma} \label{lemma_martingal}
We have $\Lambda^n\psi'=0$ for $n\geq 1$ and $\Et(\psi'\circ f^n|\Bc_m)=0$ for $m>n\geq 0$.
\end{lemma}
\proof
We deduce from the definition of $\psi''$ that
$$\Lambda\psi'=\Lambda\psi -\Lambda \psi''+\Lambda(\psi''\circ f)
=\Lambda\psi-\Lambda\psi''+\psi''=0.$$
It follows that $\Lambda^n\psi'=0$ for $n\geq 1$. 
For every function $\phi$ in $L^2(\mu)$, since $\mu$ is invariant, we have
$$\langle \mu, (\psi'\circ f^n)(\phi\circ f^m) \rangle  = 
\langle \mu, \psi'(\phi\circ f^{m-n})\rangle = \langle \mu, (\Lambda^{m-n}
\psi')\phi\rangle =0,$$
which completes the proof.
\endproof

Given a function $h\in L^1(\mu)$,  define {\it the Birkhoff's sum}
$\St_nh$ by  
$$ \St_0h:=0\qquad \mbox{and} \qquad \St_nh:= \sum\limits_{j=0}^{n-1}
h\circ f^j \quad \mbox{for } n\geq 1.$$

\begin{lemma} \label{lemma_ldt}
The coboundary $\psi''-\psi''\circ f$ satisfies the LDT.
\end{lemma}
\proof
By Proposition \ref{prop_iterate_dsh}, up to multiplying $\psi$ by a constant, we can
assume that $\langle \mu, e^{|\psi''|}\rangle \leq c$ for some
constant $c>0$. 
Observe that $\St_n(\psi''-\psi'' \circ f) =\psi''-\psi''\circ f^{n}$. 
Consequently, for  a given $\epsilon >0$, we have using the invariance
of $\mu$
\begin{eqnarray*}
\mu\big\{ | \St_n(\psi''-\psi''\circ f)|  >n\epsilon \big\}& \leq&
\mu \Big\{|  \psi''\circ f^{n}|> {n\epsilon \over 2}\Big\}+ \mu
\Big\{|\psi''|>{n\epsilon\over 2}\Big\}\\
&=& 2\mu \Big\{|  \psi''|>\frac{n\epsilon}{2}\Big\}=
2\mu \big\{ e^{| \psi''|}>e^{\frac{n\epsilon}{2}}\big\} \\
& \leq &  2e^{-\frac{n\epsilon}{2}}\langle\mu, e^{| \psi''|}\rangle \leq
2c e^{-n\epsilon\over 2}.
\end{eqnarray*}
Hence, $\psi''\circ f-\psi''$  satisfies the LDT.
\endproof

It remains to show that  $\psi'$ satisfies the LDT. Fix a number $\delta$
such that $1<\delta^5<d$. We will use the following lemma for a positive
constant $b$ of order $O(\log n)$. 

\begin{lemma} \label{lemma_unbounded}
There are constants $c>0$ and $\alpha>0$ such that
  for every $b\geq 1$ we have 
$$\mu\big\{|\psi'|>b\big\} \leq c e^{-\alpha\delta^b}.$$
\end{lemma}
\proof
Define $\varphi:=\sum_{n\geq 1} \delta^{5n} |\Lambda^n\psi|$. Since
$\|\Lambda \phi\|_\DSH'\leq d^{-1} \|\phi\|_\DSH'$ on $\DSH_0(\P^k)$,
$\varphi$ is in $\DSH_0(\P^k)$ and has a bounded DSH-norm. So,  $\langle \mu,
e^{\alpha\varphi}\rangle\leq c'$ for some constants $c'>0$ and
$\alpha>0$. 
It is enough to consider the case where $b=5l$ for some integer $l$. Since
$|\psi|\leq 1$, we also have $|\Lambda^n\psi|\leq 1$. Hence
$$|\psi''| \leq  \sum_{n\geq 1} |\Lambda^n\psi| \leq  
\delta^{-5l} \sum_{n\geq 1} \delta^{5n}| \Lambda^n\psi|
+ \sum_{1\leq n\leq l} |\Lambda^n\psi|\\
\leq \delta^{-5l} \varphi+ l.$$
Consequently, we deduce that
$$\mu\big\{|\psi''|>2l\big\}\leq \mu\big\{\varphi> \delta^{5l} \big\}\\
\leq  e^{-\alpha \delta^{5l}} \langle \mu, e^{\alpha \varphi}\rangle \leq  c'e^{-\alpha \delta^{5l}}.
$$
Therefore, by definition of $\psi'$, since $|\psi|\leq
1\leq l$ and $\mu$ is invariant, we obtain
\begin{eqnarray*}
\mu \big\{|\psi'|>5l\big\} & \leq &  \mu\big\{|\psi''|>2l\big\}
+\mu\big\{|\psi''\circ f|>2l\big\}  \\
& = & 2\mu\big\{|\psi''|>2l\big\} \leq 2 c'e^{-\alpha \delta^{5l}}.
\end{eqnarray*}
This implies the lemma.
\endproof

In order to apply Lemma \ref{lemma_Bennett}, we will need the
following property.

\begin{lemma} \label{lemma_algebra}
There is a Borel set $A$ such that $\mu(A\cap B)=(1-d^{-1})\mu(B)$ for
every $B$ in $\Bc_1$ .
\end{lemma}
\proof
Recall that $f$ defines a ramified covering of degree $d^k$.
Since $\mu$ has no mass on analytic sets, it does not charge the
critical values of $f$. So, there is a Borel set $Z$ of total $\mu$
measure such that $f^{-1}(Z)$ is the union of $d^k$ disjoint Borel
sets $Z_i$, $1\leq i\leq d^k$. Moreover, one can choose $Z_i$ so
that $f:Z_i\rightarrow Z$ is bijective. Since $\mu$ is totally
invariant, we have $\mu(Z_i)=d^{-k}$ for every $i$. Define
$A:=\cup_{i>d^{k-1}} Z_i$. 

Since $B$ is an element of $\Bc_1$, we have
$B=f^{-1}(B')$ with $B':=f(B)$. We also have 
$\mu(Z_i\cap f^{-1}(B'))=d^{-k} \mu(B')=d^{-k}\mu(B)$. Therefore,
$$\mu(A\cap B)= \sum_{i>d^{k-1}}\mu(Z_i\cap
f^{-1}(B'))=\sum_{i>d^{k-1}}d^{-k}\mu(B)=(1-d^{-1})\mu(B).$$ 
This gives the lemma.
\endproof

\begin{lemma} \label{lemma_bounded_part}
For every $b\geq 1$, there are Borel sets $W_n$ such that
$\mu(W_n)\leq c ne^{-\alpha \delta^b}$ and 
$$\int_{\P^k\setminus W_n} e^{\lambda \St_n\psi'}d\mu\leq
d \Big [{(d-1) e^{-\lambda b}+e^{(d-1)\lambda b}\over d}\Big]^n.$$
where $c>0$ and $\alpha>0$  are constants independent of $b$.
\end{lemma}
\proof
For $n=1$, define $W:=\{|\psi'|>b\}$, $W':=f(W)$ and
$W_1:=f^{-1}(W')$.
Since $\mu$ is totally invariant and $f$ has topological degree $d^k$,
we have $\mu(f(W))\leq d^k\mu(W)$. This 
and Lemma \ref{lemma_unbounded} imply that
$$\mu(W_1)=\mu(W')\leq d^k \mu(W)\leq ce^{-\alpha\delta^b}$$
for some constant $c>0$.
We also have
$$\int_{\P^k\setminus W_1} e^{\lambda \St_1\psi'}d\mu
=\int_{\P^k\setminus W_1} e^{\lambda \psi'}d\mu\leq e^{\lambda b} \leq
d \Big[{(d-1)e^{-\lambda b}+e^{\lambda b}\over d}\Big].$$
So, the lemma holds for $n=1$. 

Suppose the lemma for $n\geq 1$, we need  to prove  it for $n+1$.
Define  $W_{n+1}:=f^{-1}(W_n) \cup W_1=f^{-1}(W_n\cup W')$. 
We have
$$\mu(W_{n+1}) \leq  \mu(f^{-1}(W_n))+ \mu(W_1) = \mu(W_n)+
\mu(W_1)\leq c(n+1)e^{-\alpha\delta^b}.$$ 
We will apply Lemma \ref{lemma_Bennett} to $M:=\P^k$, $\Fc:=\Bc$, $\Gc:=\Bc_1=f^{-1}(\Bc)$,
$m:=\mu$, $\nu:=1-d^{-1}$ (see Lemma \ref{lemma_algebra}) 
and to the function $\psi^*$ such that 
$\psi^*=\psi'$ on $\P^k\setminus W_1$ and $\psi^*=0$ on $W_1$. 
By Lemma \ref{lemma_martingal}, we have
$\Et(\psi^*|\Gc)=\Et(\psi^*|\Bc_1)=0$ since $W_1$ is an element of $\Bc_1$. 

Observe that $|\psi'|\leq b$ on $\P^k\setminus W_1$. 
Hence, $|\psi^*|\leq b$. 
By Lemma \ref{lemma_Bennett}, we have 
$$\Et(e^{\lambda\psi^*}|\Bc_1)\leq  {(d-1)e^{-\lambda b}+e^{(d-1)\lambda
    b}\over d} \quad \mbox{on}\quad \P^k\quad \mbox{for}\quad \lambda\geq
0.$$ 
It follows that 
$$\Et(e^{\lambda\psi'}|\Bc_1)\leq  {(d-1)e^{-\lambda b}+e^{(d-1)\lambda
    b}\over d} \quad \mbox{on}\quad \P^k\setminus W_1\quad \mbox{for}\quad \lambda\geq
0.$$ 
Now, using the fact that $W_{n+1}$ and 
$e^{\lambda  \St_{n}(\psi'\circ f)}$ are $\Bc_1$-measurable, we can write
\begin{eqnarray*}
\int_{\P^k\setminus W_{n+1}} e^{\lambda  \St_{n+1}\psi'}d\mu
&=&\int_{\P^k\setminus W_{n+1}} 
e^{\lambda\psi'} e^{\lambda  \St_{n}(\psi'\circ f)}d\mu\\
&=&\int_{\P^k\setminus W_{n+1}} 
\Et( e^{\lambda\psi'}|\Bc_1) e^{\lambda  \St_{n}(\psi'\circ f)}d\mu.
\end{eqnarray*}
Since $W_{n+1}=f^{-1}(W_n)\cup W_1$, the last integral is bounded by
\begin{eqnarray*}
\lefteqn{\sup_{\P^k\setminus W_1} \Et( e^{\lambda\psi'}|\Bc_1)\int_{\P^k\setminus f^{-1}(W_n)}  
e^{\lambda \St_{n}(\psi'\circ f)}d\mu}\\ 
&\leq &  \Big[{(d-1)e^{-\lambda b}+e^{(d-1)\lambda b}\over d}\Big] 
\int_{\P^k\setminus W_n} 
e^{\lambda  \St_{n}\psi'}d\mu\\
&\leq &  d\Big[{(d-1)e^{-\lambda b}+e^{(d-1)\lambda b}\over d}\Big]^{n+1},
\end{eqnarray*}
where the last inequality  follows 
from the hypothesis of induction for $n$.
So, the lemma holds for $n+1$.
\endproof

The following lemma, together with Lemma \ref{lemma_ldt}, implies Theorem \ref{th_deviation}.

\begin{lemma} \label{lemma_ldt_bis}
The function $\psi'$ satisfies the LDT.
\end{lemma}
\proof
Fix an $\epsilon>0$. By Lemma \ref{lemma_bounded_part},
we  have, for  every $\lambda\geq 0$
\begin{eqnarray*}
\mu\big\{|\St_n\psi'|\geq  n\epsilon\big\}&\leq&  
\mu (W_n)+e^{-\lambda n\epsilon}
\int_{\P^k\setminus W_n} e^{\lambda  \St_n\psi'}d\mu\\
&\leq &  cne^{-\alpha \delta^b}+de^{-\lambda n\epsilon} \Big [{(d-1)
  e^{-\lambda b}
+e^{(d-1)\lambda b}\over d}\Big]^n.
\end{eqnarray*}
Take $b:=\log n (\log \delta)^{-1}$ and $\lambda:=u\epsilon b^{-2}$
with a fixed $u>0$ small enough. We have 
$$cne^{-\alpha \delta^b} = cn e^{-\alpha n} \leq e^{-\alpha n/2}$$
for $n$ large. Since $u$ is small, $\lambda b$ is small. It follows
that
\begin{eqnarray*}
{(d-1) e^{-\lambda b}+e^{(d-1)\lambda b}\over d} & \leq &  {(d-1) (1-\lambda
  b +\lambda^2b^2)+(1+(d-1)\lambda b+ (d-1)^2\lambda^2 b^2)\over d} \\
& \leq &  1+d^2\lambda^2b^2 \leq e^{d^2\lambda^2b^2} = e^{d^2u^2\epsilon^2b^{-2}}.
\end{eqnarray*}
Therefore
$$de^{-\lambda n\epsilon} \Big [{(d-1) e^{-\lambda b}+e^{(d-1)\lambda
    b}\over d}\Big]^n \leq de^{-nu\epsilon^2 b^{-2}(1-d^2u)}=
de^{-n(\log n)^{-2}h_\epsilon}$$
for some constant $h_\epsilon>0$. We deduce from the previous
estimates that
$$\mu\big\{|\St_n\psi'|\geq  n\epsilon\big\}\leq e^{-n(\log
  n)^{-2}h_\epsilon}$$
for some constant $h_\epsilon>0$ and for $n$ large. So, $\psi'$
satisfies the LDT. 
\endproof

Now, using Proposition \ref{prop_iterate_holder} we can prove the LDT for H{\"o}lder
continuous observables.

\begin{theorem} \label{th_deviation_bis}
Let $f$ be a holomorphic endomorphism of $\P^k$ of algebraic degree
$d\geq 2$. Then the equilibrium measure $\mu$ of $f$ satisfies the
large deviations theorem for H{\"o}lder continuous observables. More
precisely, if $\psi$ is a H{\"o}lder continuous function then for every
$\epsilon>0$ there is a constant $h_\epsilon>0$ such that 
\begin{equation*}
\mu\Big\{ z\in \P^k:\   \big\vert\frac{1}{n}\sum_{j=0}^{n-1}\psi\circ
f^j(z)- \langle \mu,\psi\rangle \big\vert >\epsilon \Big\}\leq
e^{-n (\log n)^{-2}h_\epsilon}
\end{equation*}
for all $n$ large enough.
\end{theorem}

The proof follows along the same lines of Theorem \ref{th_deviation}. Fix a
$\nu$-H{\"o}lder continuous function $\psi$ with $0<\nu\leq 2$ and a
constant $1<\delta<d^{\nu/10}$. We define as above the function
$\psi'$, $\psi''$ and $\varphi:=\sum_{n\geq 0} \delta^{5n}|\Lambda^n\psi|$. 
We only have to check that $\langle \mu, e^{\alpha\varphi}\rangle \leq
c$ for some constants $\alpha>0$ and $c>0$. In fact, this implies the inequality
$\langle \mu, e^{\alpha|\psi''|}\rangle \leq c$ and the crucial estimate in
Lemma \ref{lemma_unbounded}. We deduce the estimate  
$\langle \mu, e^{\alpha\varphi}\rangle \leq c$ from Proposition
\ref{prop_iterate_holder} and the following lemma for
$\theta:=\delta^5d^{-\nu/2}$ and $\eta_n:=\alpha d^{n\nu/2}|\Lambda^n\psi|$.

\begin{lemma}
Let $\eta_n$  be positive measurable functions and $0<\theta<1$ be a
constant. Assume there is a constant $c>0$ such that $\langle\mu,
e^{\eta_n}\rangle \leq c$ for every $n\geq 0$. If
$\xi:=(1-\theta)\sum_{n\geq 0} \theta^n\eta_n$, then $e^\xi$ is
$\mu$-integrable.   
\end{lemma} 
\proof
Define $\xi_m:=(1-\theta)\sum_{n\geq m} \theta^{n-m}\eta_n$. We have
$\xi_0=\xi$ and $\xi_m=(1-\theta)\eta_m+\theta\xi_{m+1}$. The H{\"o}lder
inequality implies that
$$\langle \mu, e^{\xi_m}\rangle = \langle \mu,
e^{(1-\theta)\eta_m}e^{\theta\xi_{m+1}}\rangle \leq \langle \mu,
e^{\eta_m}\rangle^{1-\theta}\langle \mu, e^{\xi_{m+1}}\rangle^\theta
\leq c^{1-\theta}\langle \mu, e^{\xi_{m+1}}\rangle^\theta.$$
By induction, this implies that
$$\langle \mu,e^{\xi_0}\rangle \leq
c^{(1-\theta)(1+\theta+\theta^2+\cdots)},$$
which implies the result.
\endproof

\begin{remark} \rm
In the main estimate of Theorem \ref{th_deviation_bis}, we can remove the factor
$(\log n)^{-2}$ if $\|\Lambda^n \psi\|_{L^\infty(\mu)}$ tends to 0
exponentially fast when $n\rightarrow \infty$. This is the case in
dimension 1, see Drasin-Okuyama
\cite{DrasinOkuyama} and when $f$ is a {\bf generic} map 
in higher dimension, see \cite{DinhSibony8}. LDT was recently proved for
Lipschitz observables in dimension 1 by Xia-Fu \cite{XiaFu}.
It seems there is a slip in their paper: they state the main result
for H\"older continuous observables. 
\end{remark}


\section{Abstract version of large deviations theorem}

In this section, we give a version of the large deviations theorem in
an abstract setting. Let $(M,\Fc,m)$ be a probability space and
$f:M\rightarrow M$ a measurable map which preserves $m$,
i.e. $f_*(m)=m$. Define $\Fc_1:=f^{-1}(\Fc)$.   We say that $f$ has {\it bounded jacobian} if
there is a constant $\kappa>0$ such that $m(f(A))\leq \kappa m(A)$
for every $A\in\Fc$. 
Observe that $f^*$ defines a linear operator of norm 1 from $L^2(m)$
into itself. 

\begin{theorem}
Let $f:(M,\Fc,m)\rightarrow (M,\Fc,m)$ be a map with bounded
jacobian which preserves $m$ as above. Let $\Lambda$ denote the adjoint of $f^*$
on $L^2(m)$. Let $\psi$ be a bounded real-valued measurable function. Assume there
are constants $\delta>1$ and $c>0$ such that 
$$\langle m,
e^{\delta^n|\Lambda^n\psi-\langle m,\psi\rangle|}\rangle \leq c\quad \mbox{for every}\quad
n\geq 0.$$ 
Then $\psi$ satisfies the large deviations theorem, that is, for every
$\epsilon>0$, there exists a constant $h_\epsilon>0$ such that 
\begin{equation*}
m\Big\{ z\in M:\   \big\vert\frac{1}{n}\sum_{j=0}^{n-1}\psi\circ
f^j(z)- \langle m,\psi\rangle \big\vert >\epsilon \Big\}\leq
e^{-n (\log n)^{-2}h_\epsilon}
\end{equation*}
for all $n$ large enough.
\end{theorem} 

The proof follows the same steps as in Section \ref{section_ldt}. The
details are left to the reader. We only notice two important points.
The property that $f$ is of bounded jacobian allows to prove an
analog of Lemma \ref{lemma_bounded_part}. Indeed, in the proof of that
lemma, the inequality
$\mu(W')\leq d^k\mu(W)$ should be replaced by $m(W')\leq\kappa
m(W)$. The 
following version of the Bennett's inequality  replaces Lemma \ref{lemma_Bennett}.

\begin{lemma}
Let $\psi$ be a real-valued function on $M$ such that
$\|\psi\|_{L^\infty(m)}\leq b$ and $\Et(\psi|\Fc_1)=0$.  
Then  
$$\Et(e^{\lambda\psi}|\Fc_1) \leq { e^{-\lambda b}+e^{\lambda b} \over 2}$$ 
for every $\lambda\geq 0$. 
\end{lemma}
\proof
We decompose the measure $m$ using the fibers of $f$. For
$m$-almost every $x\in M$, there is a
positive measure $m_x$ on $M_x:=f^{-1}(x)$ such that if $\varphi$ is a
function in $L^1(m)$ then
$$\langle m,\varphi\rangle = \int_M \langle m_x,\varphi\rangle
dm(x).$$
Since $m$ is invariant by $f$, we have
$$\langle m,\varphi\rangle = \langle m,\varphi\circ f\rangle =\int_M
\langle m_x,\varphi\circ f\rangle dm(x)=\int_M \|m_x\|\varphi(x) dm(x).$$
Therefore,  $m_x$ is a probability measure for $m$-almost every $x$.
Using also the invariance of $m$, we obtain for $\varphi$ and $\phi$ in $L^2(m)$ that
\begin{eqnarray*}
\langle m ,\varphi(\phi\circ f)\rangle & = &  \int_M \langle
m_x,\varphi(\phi\circ f)\rangle
dm(x)= \int_M \langle
m_x,\varphi\rangle \phi(x) dm(x)\\
& = & \int_M \langle
m_{f(x)},\varphi\rangle \phi(f(x)) dm(x).
\end{eqnarray*}
We deduce that
$$\Et(\varphi|\Fc_1)(x)=\langle m_{f(x)},\varphi\rangle.$$

By hypothesis, we have $\langle m_x,\psi\rangle=0$ for $m$-almost
every $x$. It suffices to check that 
$$\langle
m_x,e^{\lambda\psi}\rangle \leq { e^{-\lambda b}+e^{\lambda b} \over 2}$$ 
Consider first the particular case where there is an element
$A\subset\Fc$  such that $A\subset M_x$ and
$m_x(A)=1/2$. Applying 
Lemma \ref{lemma_Bennett} to $M_x:=f^{-1}(x)$, $m_x$, $A$, $\nu:=1/2$ and for
$\Gc:=\{\varnothing, M_x\}$ the trivial $\sigma$-algebra of $M_x$, yields the result.
The general case is deduced from the previous particular case. Indeed,
it is enough to apply this case to the disjoint union of $(M,\Fc,m)$ with a
copy $(M',\Fc',m')$ of this space, i.e. to the space $(M\cup
M',\Fc\cup\Fc',{m\over 2}+{m'\over 2})$ and to the function equal to
$\psi$ on $M$ and on $M'$.
\endproof

\small

T.-C. Dinh, UPMC Univ Paris 06, UMR 7586, Institut de
Math{\'e}matiques de Jussieu,\break F-75005 Paris, France. {\tt
  dinh@math.jussieu.fr}, {\tt http://www.math.jussieu.fr/$\sim$dinh} 

\

\noindent
V.-A.  Nguy{\^e}n, Vietnamese Academy  of Science  and  Technology,
Institute of Mathematics, Department of Analysis,
18  Hoang Quoc Viet  Road, Cau Giay  District,
10307 Hanoi, Vietnam. {\tt nvanh@math.ac.vn}

\

\noindent
N. Sibony,
Universit{\'e} Paris-Sud, Math{\'e}matique - B{\^a}timent 425, 91405
Orsay, France. {\tt nessim.sibony@math.u-psud.fr} 
\end{document}